\documentclass[12pt]{article}%

\usepackage{eurosym}
\usepackage{amsthm,amssymb,amsmath}
\usepackage{epsfig}
\usepackage{xcolor}
\usepackage{pgf}
\usepackage{tikz,fullpage}
\usepackage{tkz-berge}
\usepackage[latin1]{inputenc}
\usepackage{amsmath}
\usepackage{amsfonts}
\usepackage{amssymb}
\usepackage{graphicx}%
\newtheorem{prelem}{{\bf Theorem}}

\newtheorem{theorem}{Theorem}
\newtheorem{corollary}[theorem]{Corollary}
\newtheorem{lemma}[theorem]{Lemma}
\newtheorem{obs}[theorem]{Observation}
\newtheorem{prop}[theorem]{Proposition}

\theoremstyle{definition}

\theoremstyle{remark}
\newtheorem{rem}[theorem]{Remark}
\DeclareMathOperator{\opt}{\emph{Opt}}
\begin{document}

\title{On the outer independent total double Roman dominating functions}
\date{}
\author{$^{(1)}$H.A. Ahangar, $^{(2)}$M. Chellali,\\$^{(3)}$S.M. Sheikholeslami and $^{(4)}$J.C. Valenzuela-Tripodoro\vspace
{7.5mm} \\$^{(1)}${\small Department of Basic Sciences}\\
{\small \ Babol University of Technology}\\
{\small \ Babol, Iran.}\\
{\small \ ha.ahangar21@gmail.com\vspace{2mm}}\\
$^{(2)}${\small Department of Basic Sciences}\\
{\small \ University of Blida}\\
{\small \ Blida, Algeria}\\
{\small \ m\_chellali@yahoo.com\vspace{2mm}}\\
$^{(3)}${\small Department of Basic Sciences}\\
{\small \ Azarbaijan Shahid Madani University }\\
{\small \ Tabriz, Iran}\\
{\small \ sm\_sheikholeslami@yahoo.com\vspace{3mm}}\\
$^{(4)}${\small Department of Basic Sciences, University of C\'{a}diz, Spain.}
\\
\ {\small jcvtrip@gmail.com\vspace{5mm}}}
\maketitle

\begin{abstract}
Let $\{0,1,\dots, t\}$ be abbreviated by $[t].$ A double Roman dominating function (DRDF) on a graph $\Gamma=(V,E)$ is a map
$l:V\rightarrow [3]$ satisfying \textrm{(i)} if $l(r)=0$ then there
must be at least two neighbors labeled 2 under $l$ or a neighbor $r'$ with
$l(r')=3$; and \textrm{(ii)} if $l(r)=1$ then $r$ must be adjacent to a vertex
$r'$ such that $l(r')\geq2$. A DRDF is an outer-independent total double Roman
dominating function (OITDRDF) on $\Gamma$ if the set of vertices labeled $0$
induces an edgeless subgraph and the subgraph induced by the vertices with a
non-zero label has no isolated vertices. The weight of an OITDRDF is the sum
of its map values over all vertices, and the outer independent total
Roman dominating number $\gamma_{tdR}^{oi}(\Gamma)$ is the minimum weight of an
OITDRDF on $\Gamma$. First, we prove that the problem of determining $\gamma
_{tdR}^{oi}(\Gamma)$ is NP-complete for bipartite and chordal graphs, after that, we prove
that it is solvable in linear time when we are restricting to bounded
clique-width graphs. Moreover, we present some tight bounds on $\gamma
_{tdR}^{oi}(\Gamma)$ as well as the exact values for several graph families.


\end{abstract}

\section{Introduction}

The starting point of Roman domination in graphs, as well as all its variants,
can be attributed to the mathematical formalization of a defensive model of
the Roman Empire described by Stewart in \cite{Stewart}. The formal definition
was given by Cockayne et al. \cite{CDHH} as follows. Given a graph $\Gamma=(V,E),$
with vertex set $V$ and edge set $E$, a \emph{Roman dominating function} (RDF)
on $\Gamma$ is a map $l:V\rightarrow [2]$, that labels labels to
vertices of $\Gamma$, such that each vertex labeled with $0$ must be adjacent to a
vertex with a label $2$. The sum of all vertex labels, $w(l)=l(V)=\sum_{r\in
V}l(r)$, is called the weight of the RDF $l$ and the minimum weight over all
possible RDF's is the \emph{Roman domination number}, $\gamma_{R}(\Gamma)$, of the
graph $\Gamma$. For the sake of simplicity, {an }RDF with minimum weight is known
as a $\gamma_{R}(\Gamma)$-function {or a }$\gamma_{R}$-function of $\Gamma$. Evidently,
there is a one-on-one relation between Roman dominating functions and the set
of subsets $\{V_{0}^{l},V_{1}^{l},V_{2}^{l}\}$ of $V(\Gamma)$, where $V_{i}%
^{l}=\{r\in V|l(r)=i\}.$ That is why {an} RDF $l$ is usually represented as
$l=(V_{0}^{l},V_{1}^{l},V_{2}^{l})$ or simply by $(V_{0},V_{1},V_{2})$, if
there is no possibility of confusion.

Following the publication of the introductory article of Roman domination,
more than 200 articles have been published, among which are devoted to the
study of interesting variations of the original problem. For example, we can
highlight \emph{mixed Roman domination} \cite{ahv}, \emph{double Roman
domination} \cite{bhh}, \emph{strong Roman domination} \cite{amsvy},\emph{
total Roman domination} \cite{acs2}, \emph{outer-independent Roman domination}
\cite{ahsy}, \emph{outer-independent total Roman domination} \cite{cky} and
\emph{outer-independent double Roman domination} \cite{acs3}. { The
reader can consult the following book chapters \cite{Ch1, Ch2} and the survey
\cite{Ch3} to delve into the main results on Roman domination and its
variants.}

{ We say that a map
$l:V(\Gamma)\rightarrow\lbrack3]$ is a  DRDF,
 if it meets the requirements described below: }

\begin{itemize}
\item Any vertex $r\in V$ with $l(r)=0$ is adjacent to either a vertex $y$
such that $l(y)=3$ or to two vertices $y,y^{\prime}\in V$ with
$l(y)=l(y^{\prime})=2.$

\item Any vertex $r\in V$ with $l(r)=1$ is adjacent to a vertex $y\in V$
with $l(y)\ge2.$
\end{itemize}

The \textit{double Roman domination number} (\emph{DRD-number for short} )
$\gamma_{dR}(\Gamma)$ equals the minimum weight of a DRDF on $\Gamma$, and a DRDF of $\Gamma$
with weight $\gamma_{dR}(\Gamma)$ is called a $\gamma_{dR}(\Gamma)$-function$.$ A
\emph{total double Roman dominating function} (TDRDF) is a DRDF $(V_{0}%
,V_{1},V_{2},V_{3})$ for which the subgraph induced by $V-V_{0}$ has no
isolated vertices. Likewise, an  OIDRDF is a DRDF $(V_{0},V_{1},V_{2},V_{3})$ for which
$V_{0}$ induces an edgeless subgraph. Let $\gamma_{tdR}(\Gamma)$ and $\gamma
_{oidR}(\Gamma)$ denote the \textit{total double Roman domination number}
(\emph{TDRD-number} for short) and the \textit{outer-independent Roman
domination number }(\emph{OIDRD-number} for short) of $\Gamma,$ respectively.

Our aim in this work is to start the study of a new variation of Roman
domination, namely the Outer independent total double Roman domination. In
this case, a DRDF $l=(V_{0},V_{1},V_{2},V_{3})$ is an OITDRDF if $V_{0}$ induces an
edgeless subgraph and $V_{1}\cup V_{2}\cup V_{3}$ induce an isolated-free
vertex subgraph. Analogously, we define the \emph{outer independent total
double Roman domination number }(\emph{OITDRD-number} for short),
$\gamma_{tdR}^{oi}(\Gamma)$, as the minimum weight of an OITDRDF on $\Gamma$. Every
OITDRDF on $\Gamma$ having minimum weight is called a $\gamma_{tdR}^{oi}(\Gamma)$-function.

{  As a matter of example, one can model a sewerage network by
means of a graph. Sewage consists of the used water that comes from
commercial, domestic, or industrial sources which contain mostly water and
less than $1\%$ of dissolved or suspended solids. In the underlying graph, we
represent each household or business by a vertex with a label $0$. The network
of pipes (edges of the graph) takes the sewage away from the $0$-vertices to
the sewers. We can have two different classes of sewers: collector sewers,
with small waste collection capacity and main sewers, with a higher capacity.
Collector (resp. main) sewers are labeled with $2$ (resp. $3$) as vertices of
the graph. Every individual household or business must be connected to either
two collector sewers or either to a main sewer. Most sewerage systems are
designed to use gravity but there are many situations in which it is not
possible. So, we can find sewage pumps in the network that should be directly
connected to a collector/main sewer. We represent a sewage pump by a
$1$-vertex in the graph. The service branch lines from individual sewage
sources are disjoint to avoid interferences, that is, no two different
$0$-vertices are connected by an edge. Moreover, to guarantee backup against
overload in operation, each collector/main sewer must be connected to, at
least, one other collector/main sewer of the same class or of a different one.
In graph terminology, the latter means that the set of vertices with a label
greater than or equal to $1$ induce a subgraph with no isolated vertex.
Evidently, a sewage network satisfies these conditions iff the corresponding
labeling of the underlying graph is an OITDRDF. So, we can minimize the cost of the network by
determining the OITDRD-number of the graph. }

\begin{figure}[ptb]
\centering
\resizebox{\textwidth}{!}{
\begin{tikzpicture}[scale=1,transform shape,font=\Large]
\xdefinecolor{borde}{RGB}{0,0,0}
\xdefinecolor{texto}{RGB}{0,0,0}
\xdefinecolor{relleno}{RGB}{255,255,255}
\tikzstyle{VertexStyle}=[shape=circle,color=texto!100,draw=borde!100,minimum size=18mm,fill=relleno!100]
\Vertex[x=120mm,y=340mm,LabelOut=false,L=\hbox{H1}]{0}
\tikzstyle{VertexStyle}=[shape=circle,color=texto!100,draw=borde!100,minimum size=18mm,fill=relleno!100]
\Vertex[x=200mm,y=360mm,LabelOut=false,L=\hbox{H2}]{1}
\tikzstyle{VertexStyle}=[shape=circle,color=texto!100,draw=borde!100,minimum size=18mm,fill=relleno!100]
\Vertex[x=280mm,y=340mm,LabelOut=false,L=\hbox{H3}]{2}
\tikzstyle{VertexStyle}=[shape=circle,color=texto!100,draw=borde!100,minimum size=18mm,fill=relleno!100]
\Vertex[x=200mm,y=200mm,LabelOut=false,L=\hbox{H4}]{3}
\tikzstyle{VertexStyle}=[shape=circle,color=texto!100,draw=borde!100,minimum size=18mm,fill=relleno!100]
\Vertex[x=120mm,y=220mm,LabelOut=false,L=\hbox{H5}]{4}
\tikzstyle{VertexStyle}=[shape=circle,color=texto!100,draw=borde!100,minimum size=18mm,fill=relleno!100]
\Vertex[x=280mm,y=220mm,LabelOut=false,L=\hbox{H6}]{5}
\tikzstyle{VertexStyle}=[shape=circle,color=texto!100,draw=borde!100,minimum size=29mm,fill=relleno!100]
\Vertex[x=200mm,y=280mm,LabelOut=false,L=\hbox{Sewage pump}]{6}
\tikzstyle{VertexStyle}=[shape=rectangle,color=texto!100,draw=borde!100,minimum size=30mm,fill=relleno!100]
\Vertex[x=160mm,y=300mm,LabelOut=false,L=\hbox{Collector sewer}]{7}
\tikzstyle{VertexStyle}=[shape=rectangle,color=texto!100,draw=borde!100,minimum size=30mm,fill=relleno!100]
\Vertex[x=240mm,y=300mm,LabelOut=false,L=\hbox{Collector sewer}]{8}
\tikzstyle{VertexStyle}=[shape=rectangle,color=texto!100,draw=borde!100,minimum size=30mm,fill=relleno!100]
\Vertex[x=200mm,y=240mm,LabelOut=false,L=\hbox{Main sewer}]{9}
\xdefinecolor{pincel}{RGB}{0,0,0}
\Edge[label=,style={color=pincel, ,line width=1.7mm}](6)(8)
\Edge[label=,style={color=pincel, ,line width=1.7mm}](7)(6)
\Edge[label=,style={color=pincel, ,line width=1.7mm}](9)(6)
\Edge[label=,style={color=pincel, ,line width=1.7mm}](1)(7)
\Edge[label=,style={color=pincel, ,line width=1.7mm}](0)(7)
\Edge[label=,style={color=pincel, ,line width=1.7mm}](2)(7)
\Edge[label=,style={color=pincel, ,line width=1.7mm}](2)(8)
\Edge[label=,style={color=pincel, ,line width=1.7mm}](1)(8)
\Edge[label=,style={color=pincel, ,line width=1.7mm}](0)(8)
\Edge[label=,style={color=pincel, ,line width=1.7mm}](9)(4)
\Edge[label=,style={color=pincel, ,line width=1.7mm}](9)(3)
\Edge[label=,style={color=pincel, ,line width=1.7mm}](9)(5)
\tikzstyle{VertexStyle}=[shape=circle,color=texto!100,draw=borde!100,minimum size=10mm,fill=relleno!100]
\Vertex[x=320mm, y=280mm,LabelOut=false,L=\hbox{0}]{g1}
\tikzstyle{VertexStyle}=[shape=circle,color=texto!100,draw=borde!100,minimum size=10mm,fill=relleno!100]
\Vertex[x=340mm, y=320mm,LabelOut=false,L=\hbox{0}]{g2}
\tikzstyle{VertexStyle}=[shape=circle,color=texto!100,draw=borde!100,minimum size=10mm,fill=relleno!100]
\Vertex[x=340mm, y=240mm,LabelOut=false,L=\hbox{0}]{g3}
\tikzstyle{VertexStyle}=[shape=circle,color=texto!100,draw=borde!100,minimum size=10mm,fill=relleno!100]
\Vertex[x=380mm, y=280mm,LabelOut=false,L=\hbox{3}]{g4}
\tikzstyle{VertexStyle}=[shape=circle,color=texto!100,draw=borde!100,minimum size=10mm,fill=relleno!100]
\Vertex[x=420mm, y=280mm,LabelOut=false,L=\hbox{2}]{g5}
\tikzstyle{VertexStyle}=[shape=circle,color=texto!100,draw=borde!100,minimum size=10mm,fill=relleno!100]
\Vertex[x=500mm, y=280mm,LabelOut=false,L=\hbox{2}]{g6}
\tikzstyle{VertexStyle}=[shape=circle,color=texto!100,draw=borde!100,minimum size=10mm,fill=relleno!100]
\Vertex[x=460mm, y=320mm,LabelOut=false,L=\hbox{0}]{g7}
\tikzstyle{VertexStyle}=[shape=circle,color=texto!100,draw=borde!100,minimum size=10mm,fill=relleno!100]
\Vertex[x=460mm, y=240mm,LabelOut=false,L=\hbox{0}]{g8}
\tikzstyle{VertexStyle}=[shape=circle,color=texto!100,draw=borde!100,minimum size=10mm,fill=relleno!100]
\Vertex[x=460mm, y=280mm,LabelOut=false,L=\hbox{0}]{g9}
\tikzstyle{VertexStyle}=[shape=circle,color=texto!100,draw=borde!100,minimum size=10mm,fill=relleno!100]
\Vertex[x=400mm, y=240mm,LabelOut=false,L=\hbox{1}]{g10}
\Edge[label=,style={color=pincel, ,line width=1.1mm}](g1)(g4)
\Edge[label=,style={color=pincel, ,line width=1.1mm}](g2)(g4)
\Edge[label=,style={color=pincel, ,line width=1.1mm}](g3)(g4)
\Edge[label=,style={color=pincel, ,line width=1.1mm}](g10)(g4)
\Edge[label=,style={color=pincel, ,line width=1.1mm}](g10)(g5)
\Edge[label=,style={color=pincel, ,line width=1.1mm}](g10)(g6)
\Edge[label=,style={color=pincel, ,line width=1.1mm}](g4)(g5)
\Edge[label=,style={color=pincel, ,line width=1.1mm}](g5)(g7)
\Edge[label=,style={color=pincel, ,line width=1.1mm}](g5)(g8)
\Edge[label=,style={color=pincel, ,line width=1.1mm}](g5)(g9)
\Edge[label=,style={color=pincel, ,line width=1.1mm}](g6)(g7)
\Edge[label=,style={color=pincel, ,line width=1.1mm}](g6)(g8)
\Edge[label=,style={color=pincel, ,line width=1.1mm}](g6)(g9)
\end{tikzpicture}
}
\end{figure}

In what follows, we only discuss simple graphs $\Gamma=(V,E)$. The size of a graph is its number of edges and
its order is the number of elements in $V$. The \emph{open} (resp.
\emph{closed}) \emph{neighborhood} $N(r)$ (resp. $N[r]$) of a vertex $r$ is
the set $\{r'\in V(\Gamma):rr'\in E(\Gamma)\}$ (resp. $N[r]=N(r)\cup\{r\}$). The number of
adjacent vertices with $r$ is its \emph{degree}, $d_{\Gamma}(r)=|N(r)|$.
{ We denote by $\delta=\delta(\Gamma)$\text{ (resp., $\Delta=\Delta(\Gamma)$%
)} the \emph{minimum } (resp., \emph{maximum}) \emph{degree }of a graph $\Gamma.$}
A \emph{leaf} in a graph is a vertex whose degree is equal to $1$ and its
neighbor is a \emph{stem}. A stem is \emph{ strong }if it is adjacent to at
least two leaves and it is \emph{weak} otherwise. An \emph{$m$-regular graph}
is a graph such that every vertex has the same degree $m$.

We denote by $P_{p}$ the \emph{path graph} of order $p$, by $C_{p}$ the
\emph{cycle graph} of length $p$ and by $\overline{K_{p}}$ for the
\emph{edgeless graph} with $p$ vertices. { A connected graph
containing no cycle is called a \emph{tree}. A tree containing just two
vertices that are not leaves is a \emph{double star}. A double star, denoted
by $DS_{r,t},$ has $r$ and $t$ leaves attached at each non-leaf
vertex, respectively.} The \emph{corona }of a graph $\Gamma$ denoted \emph{Cor}$(\Gamma),$ is the
graph formed {from a copy} of $\Gamma$ by adding a new vertex
$y^{\prime}$ and the edge $yy^{\prime}$, for any $y\in V$.

A set of vertices $M\subseteq V$ is a \textit{dominating set} (DS) in a graph $\Gamma$
if each vertex of $V\smallsetminus M$ is adjacent to a vertex in $M.$ The
\textit{domination number} $\gamma(\Gamma)$ of a graph $\Gamma$ is the minimum
cardinality of a DS in the graph. A DS $M$ is a
\emph{total dominating set} (TDS) if any vertex in $\Gamma$ is adjacent to a vertex in
$M$ and, similarly, the minimum cardinality of a TDS is
denoted by $\gamma_{t}(\Gamma)$. A TDS $M$ of a graph $\Gamma$ is a
\emph{total co-independent dominating set (Tco-IDS)}if the set $V\smallsetminus M$ is
non-empty and induces an edgeless subgraph. As usual, the minimum cardinality
of such a Tco-IDS is denoted by $\gamma
_{t,coi}(\Gamma).$ Evidently for each graph $\Gamma$ with minimum degree grater than or equal to one, we have
that
\[
\gamma(\Gamma)\leq\gamma_{t}(\Gamma)\leq\gamma_{t,coi}(\Gamma).
\]

We make use of the following results in this article.

\begin{obs}
\label{ob}\emph{If $y$ and $y'$ are stem and a leaf of a graph $\Gamma$, with $yy'\in E(\Gamma)$, then
for any OITDRDF of $\Gamma$, $s(y)+s(y')\geq3$ and $s(y)\geq1$.}
\end{obs}

\begin{theorem}
[\cite{Hao}]\label{cycle} For $p\geq3$, $\gamma_{tdR}(C_{p})=\lceil\frac
{6p}{5}\rceil$.
\end{theorem}

\begin{theorem}
[\cite{Hao}]\label{path} For $p\geq3,$ $\gamma_{tdR}(P_{p})=\left\{
\begin{array}
[c]{ll}%
6 & \mathrm{if}\;p=4\\[.5 em]%
\lceil\frac{6p}{5}\rceil & \mathrm{o.w.}.
\end{array}
\right.  $
\end{theorem}

\begin{theorem}
[\cite{Shao}]\label{tree1} \emph{For each connected graph $\Gamma$ of order $p\geq
2$, $\gamma_{tdR}(\Gamma)\leq\frac{3p}{2}$ .}
\end{theorem}

We close this part by giving the exact values of the OITDRD-number for
paths and cycles.

\begin{lemma}
\label{pathandcycle} For $p\geq3$, \newline

$(i)$ $\gamma_{tdR}^{oi}(P_{p})=\left\{
\begin{array}
[c]{lll}%
6 & \mathrm{if} & p=4,\\[.5em]%
\lceil\frac{6p}{5}\rceil &  & \mathrm{otherwise}.
\end{array}
\right.  $ \vspace{.5cm} $(ii)$ ${ \gamma_{tdR}^{oi}}(C_{p})=\lceil\frac
{6p}{5}\rceil.$
\end{lemma}

\textbf{Proof. } (i) Assume that $P_{p}=b_{1}b_{2}\ldots b_{p}$ be a path on
$p$ vertices. For $p\leq7$ the result is trivial. Therefore, we assume that{
$p\geq8$}. discuss the map $t:V(P_{p})\rightarrow\lbrack3]$ given by{
$t(b_{5j+1})=t(b_{5j+5})=1$, $t(b_{5j+2})=t(b_{5j+4})=2$ and $t(b_{5j+3})=0$
for $0\leq j\leq\left\lfloor \frac{p-5}{5}\right\rfloor $, and also
$t(b_{p})=2$, if $p\equiv1\pmod 5$; $t(b_{p})=2,t(b_{p-1})=1$ if
$p\equiv2\pmod 5;$ $t(b_{p})=t(b_{p-2})=1,t(b_{p-1})=2$ if $p\equiv3\pmod 5$;
and finally $t(b_{p})=1,t(b_{p-1})=t(b_{p-3})=2$, $t(b_{p-2})=0$, if
$p\equiv4\pmod 5$}. Evidently, $t$ is an OITDRDF of $P_{p}$ with
$\omega(t)=\lceil6p/5\rceil$, and so {$\gamma_{tdR}^{oi}(P_{p})\leq
\lceil6p/5\rceil$.}

The {equality }follows from {Theorem \ref{path}} and {the fact that
$\gamma_{tdR}^{oi}(P_{p})\geq\gamma_{tdR}(P_{p})$}.

(ii) Applying the map $t$ defined for paths, we can see that {$\gamma
_{tdR}^{oi}(C_{p})\leq\lceil6p/5\rceil$.} The {equality} follows from {Theorem
\ref{cycle} and the fact that $\gamma_{tdR}^{oi}(C_{p})\geq\gamma_{tdR}%
(C_{p})$. Thus $\gamma_{tdR}^{oi}(C_{p})=\lceil6p/5\rceil$.} $\hfill\Box$



\section{Complexity results}

The {purpose} of this {part }is to {investigate }the complexity of the
following decision problem, to which we shall refer as OITD-ROM-DOM:

\bigskip

\textbf{OITD-ROM-DOM}

\textbf{Instance}: Graph $\Gamma=(V,E)$, integer $k>0$.

\textbf{Question}: Does $\Gamma$ have an OITDRF of weight at most $k$?

\bigskip

{To prove} that this decision problem is NP-complete, we {perform }a
polynomial time reduction from the OITDRD
problem shown to be NP-complete for bipartite and chordal graphs in
\cite{acs3}. \newline\medskip

\textbf{OIDRD}

\textbf{INSTANCE:} A graph $\Gamma$ and a  integer $k>0$.

\textbf{QUESTION:} { Is $\gamma_{oidR}(\Gamma)\leq k$ }?\newline

\begin{theorem}
Problem OITD-ROM-DOM is NP-Complete for bipartite and chordal graphs\textbf{.}
\end{theorem}

\textbf{Proof. }{Evidently}, OITD-ROM-DOM is a member of $\mathcal{NP}$ since
for a given map $l=(V_{0},V_{1},V_{2},V_{3})$ of a graph $\Gamma$ we {may
verify }in polynomial time {whether }$l$ is a TDRDF for $\Gamma$ with weight at most
$k$ and that $V_{0}$ is an independent set (IS).

Now, given a  integer $k>0$ and a graph $\Gamma$ of order $p$, we construct a
graph${\ H}$ by adding for {any} vertex $x_{i}$ a star $K_{1,4},$ with center
vertex $u_{i}$ and leaves $a_{i},b_{i},c_{i},d_{i},$ attached by $a_{i}x_{i}$
at $x_{i}.$ Evidently, $\left\vert V({H})\right\vert =6\left\vert
V(\Gamma)\right\vert $ and $\left\vert E({H})\right\vert =\left\vert
E(\Gamma)\right\vert +5\left\vert V(\Gamma)\right\vert ,$ and so ${H}$ may be
constructed from $\Gamma$ in polynomial time. In addition, it is obvious that if $\Gamma$ is
a bipartite or chordal, then ${H}$ is also bipartite or chordal, respectively.

In the next, we will {prove} {that }$\Gamma$ has an OIDRDF $t$ with $\omega(t)\leq
k$ {iff }${H}$ has an OITDRDF $l$ with $\omega(l)\leq k+4p$. Suppose first
that $t$ is an OIDRDF {for} $\Gamma$ with $\omega(t)\leq k.$ Define the map
$l$ on $V({H})$ by $l(z)=t(z)$ for all $z\in V(\Gamma),$ and for {any}
$i\in\{1,2,...,p\}$, let $l(a_{1})=1,$ $l(u_{i})=3,$ $l(b_{i})=l(c_{i}%
)=l(d_{i})=0.$ {Of course}, $\omega(l)=\omega(t)+4p$ and thus $\omega(l)\leq
k+4p.$ In addition, it is obvious to notice that $l$ is an OITDRDF of ${H}.$

Conversely, suppose $l$ is an OITDRDF of ${H}$ with $\omega(l)\leq k+4p$. We
first note that for {each} $i,$ the total values labeled to $u_{i}$ and its
neighbors under $l$ is at least $4$ with $l(u_{i})=3.$ We also note that if
$l(a_{i})=0$ for some $i,$ then $l(x_{i})\neq0$ (since $l$ is an OITDRDF) and
one of $b_{i},c_{i},d_{i},$ say $b_{i},$ is labeled a non-zero value. In that
case, we can interchange the values of $a_{i}$ and $b_{i}$ and evidently $l$
remains an OITDRDF. Therefore, we can suppose w.l.g., that
$l(a_{i})>0$ for any $i.$ We claim that each $a_{i}$ can be labeled the
value $1.$ Suppose to the contrary that $l(a_{i})\in\{2,3\}$ for some $i.$ If
$l(a_{i})=3,$ then we must have $l(x_{i})=0$ and thus by relabeling $x_{i}$
and $a_{i}$ the values $2$ and $1$ respectively, the map $l$ is still an
OITDRDF of $\Gamma$ with weight at most $k+4p.$ If $l(a_{i})=2,$ then by definition
$x_{i}$ has only one neighbor in $\Gamma$ labeled a $2,$ and this case relabeling
a $1$ to $x_{i}$ and $a_{i}$ does not affect the weight of $l$ which remains
an OITDRDF of $\Gamma.$ Consequently, let us assume that $l(a_{i})=1$ for each $i.$
Now evidently, if $t$ denotes the restriction of $l$ on $V(\Gamma),$ then $t$ is an
OIDRDF of $\Gamma$ and thus $\omega(t)\leq\omega(l)-4p\leq$ $k+4p-4p=k.$ This
finishes the prove. $\ \ \hfill\Box$

\bigskip

Our next aim is to prove that this problem can be solved in linear time for
many classes of graphs. To do that, we make apply of several concepts and tools
related to monadic second-order logic (MSOL). Particularly, we use an
extension of MSOL, known as LinEMSOL, which allows to search for sets of
vertices that are optimal with respect to some evaluation map.

According to the notations and definitions described by Courcelle et
al.\cite{CMR} and Liedloff et al.\cite{LKLP}, a $m$-expression for a given
graph $\Gamma$ is an algebraic expression built applying a set of operations of the
following types: (i) $\bullet j({z})$, to add a new vertex $z$ with a label $j
\in\{1,2,\ldots,m\}$; (ii) $\Gamma_{1} \oplus G_{2}$, that is the disjoint union
of both graphs $\Gamma_{i}$; (iii) $\eta_{i,j}(\Gamma)$, that consists in the graph
obtained by adding and edge between all $i$ labeled vertices and all $j$
labeled vertices in $\Gamma$; and (iv) $\rho_{i\rightarrow j}(\Gamma)$, which is the
graph obtained from $\Gamma$ by converting every $i$ label into $j$. For example,
the complete bipartite graph $K_{3,3}$ can be described by means of the
following $2$-expression
\[
\eta_{1,2} \left(  \bullet1({a}) \oplus\bullet1({b}) \oplus\bullet1({c})
\oplus\bullet2({d}) \oplus\bullet2({e}) \oplus\bullet2({f}) \right)
\]

The minimum value $k$ such that there exists a $m$-expression of a given graph
$\Gamma$ is called the \emph{clique-width} of the graph $\Gamma$ and it is denoted by
$cwd(\Gamma).$ The clique-width of a graph is well-defined because it is always
possible to describe any graph by applying a $k$-expression.

Let us denote by $\Gamma(\tau_{1})$ the graph $\Gamma$ as a logic structure
$<V(\Gamma),\mathcal{R}>$. With this notation, we assume that $v\ \mathcal{R}\ w$,
or simply $\mathcal{R}(v,w),$ iff $v$ and $w$ are neighbors in $\Gamma$. Moreover,
let us denote by $MSOL(\tau_{1})$ the monadic second order logic with
quantifications over subsets of elements of $\Gamma(\tau_{1})$. Every optimization
problem that can be described in the following way
\begin{equation}
\opt\;\left\{  \sum_{1\leq i\leq l}a_{i}|X_{i}|\;:\;<G(\tau_{1}),X_{1}%
,\ldots,X_{l}>\;\vDash\theta(X_{1},\ldots,X_{l})\right\}  \label{optpro}%
\end{equation}
where $X_{j}$ are free variables, $a_{i}$ are integers and $\theta$ is a
formula belonging to $MSOL(\tau_{1})$, is a problem belonging to the class
$LinEMSOL(\tau_{1})$ (see~\cite{LKLP}).

We use the Theorem 30 by Liedloff et al.\cite{LKLP}, which is a version of
Theorem 4 by Courcelle et al.\cite{CMR}, stating that \emph{each LinEMSOL
optimization problem on the class of bounded clique-width graphs, $cwd(\Gamma)\le
m$, can be solved in linear time, if a $m$-expression of the graph is part of
the input}. The class of bounded clique-width graphs contains several
well-known families of graphs as, among others, distance-hereditary graphs,
series-parallel graphs or cographs. It is worth noting that each bounded
treewidth graph is a bounded clique-width graph, so all trees are also bounded
clique-width graphs.

In the rest of this part, we apply this result to show that the decision
problem associated to the OITDRDF can
be solved in linear time for the class of bounded clique-width graphs.

\begin{theorem}
The decision problem associated to the OITD-ROM-DOM problem belongs to the
class $LinEMSOL(\tau_{1}).$
\end{theorem}

\textbf{Proof. } It suffice to show that we can describe the OITD-ROM-DOM
problem as an optimization problem with the structure described in
(\ref{optpro}).

Let $l=(V_{0},V_{1},V_{2},V_{3})$ be an OITDRDF on $\Gamma$ and define the
following set of free variables $\{X_{j}\}$ where $X_{j}(v)=1$ if $v\in V_{j}$
and $X_{j}(v)=0,$ otherwise. Then, it is clear that $|X_{j}|=\sum_{v\in
V}X_{j}(v)=|V_{j}|.$ Moreover, to find an OITDRDF with minimum weight is
equivalent to solve the optimization problem,
\[
\min_{X_{i}}\left\{  |X_{1}|+2|X_{2}|+3|X_{3}|\;:\;<G(\tau_{1}),X_{0}%
,\ldots,X_{3}>\;\vDash\theta(X_{0},\ldots,X_{3})\right\}  ,
\]
where $\theta(X_{0},\ldots,X_{3})=\theta_{1}(X_{0},\ldots,X_{3})\wedge
\theta_{2}(X_{0},\ldots,X_{3})\wedge\theta_{3}(X_{0},\ldots,X_{3})$ and the
auxiliary formulae $\theta_{i}$ are defined as follows
\[%
\begin{array}
[c]{rcl}%
\theta_{1}(X_{0},\ldots,X_{3}) & = & \forall q\left[
{\ \vphantom{\frac{\frac{A}{B}}{\frac{A}{B}}}}X_{3}(q)\vee X_{2}(q)\vee\left[
{\vphantom{\frac{\int A}{B}}}X_{1}(q)\wedge\exists r\left(  \left(
X_{2}(r)\vee X_{3}(r)\right)  \wedge\mathcal{R}(q,x)\right)
{\vphantom{\frac{\int A}{B}}}\right]  \vee\right. \\[0.5em]%
\multicolumn{3}{r}{\vee\left.  \left[  {\vphantom{\frac{\int A}{B}}}%
X_{0}(q)\wedge\left(  {\vphantom{\frac{\int A}{B}}}\exists r\left(
X_{3}(r)\wedge\mathcal{R}(q,r)\right)  \vee\exists r,v\left(  X_{2}(r)\wedge
X_{2}(v)\wedge\mathcal{R}(q,r)\wedge\mathcal{R}(q,v){\vphantom{\frac{A}{B}}}%
\right)  {\vphantom{\frac{\int A}{B}}}\right)  {\vphantom{\frac{\int A}{B}}}%
\right]  {\ \vphantom{\frac{\frac{A}{B}}{\frac{A}{B} }}}\right] }\\[1.25em]%
\theta_{2}(X_{0},\ldots,X_{3}) & = & \forall q\left[  X(q)\vee\left(
{\vphantom{\frac{{\vphantom{\frac{A}{B}}}}{B}}}\exists v\left(  X_{1}(v)\vee
X_{2}(v)\vee X_{3}(v)\right)  \wedge R(q,v)\right)
{\vphantom{\frac{{\vphantom{\frac{A}{B}}}}{B}}}\right] \\[1.25em]%
\theta_{3}(X_{0},\ldots,X_{3}) & = & \forall q,r\left[
{\vphantom{\frac{{\vphantom{\frac{A}{B}}}}{B}}}X_{0}(q)\wedge X_{0}%
(r)\rightarrow\lnot R(q,r)\right]
\end{array}
\]

It is not difficult to verify that the so-defined { formula} $\theta$ checks
for the map $l$ to be an OITDRDF. The main auxiliary { formula},
$\theta_{1}$, has four clauses to check that the map satisfies the
definition of DRDF. The { formula} $\theta_{2}$ checks that $l$ is a TRDF and,
finally, $\theta_{3}$ verifies that $V_{0}$ is an IS of vertices.

Hence, $l$ is an OITDRDF iff the expression $\theta$ is assured, which
concludes the proof. $\vphantom{a}\hfill\Box$

\begin{corollary}
Problem OITD-ROM-DOM can be solved in linear time for the class of $m$-bounded
clique-width graphs whenever either there exists a linear-time algorithm to
make a $m$-presentation of the given graph or a $m$-presentation is part of
the input.
\end{corollary}

As bounded tree-width graphs are also bounded clique-width graphs, the former
corollary can be applied to this family of graphs.


\section{Bounds}

Recall that a \emph{matching} in a graph $\Gamma$ is a subset of pairwise
non-adjacent edges. The \emph{matching number} $\alpha^{\prime}(\Gamma)$
($\alpha^{\prime}$ for short) is the size of a largest matching in $\Gamma$.

\begin{theorem}
\label{Th1}\emph{Let $\Gamma$ be a connected graph of order $p\geq2$. Then the
following occurs. }

\begin{enumerate}
\item \emph{$\gamma_{tdR}^{oi}(\Gamma)\leq p+\alpha^{\prime}(\Gamma)$. Moreover, this
bound is sharp for every graph $\Gamma\in\{Cor(F)\mid F \;\text{is a connected
graph} \}$. }

\item \emph{If $\Gamma$ is triangle-free graph with $\delta(\Gamma)\geq2$, then ${
\gamma_{tdR}^{oi}(\Gamma)} \leq3\alpha^{\prime}(\Gamma)$. Moreover, this bound is tight
for every graph $\Gamma\in\{ Cor(F)\mid F \; \text{is a connected triangle-free
graph} \}$. }
\end{enumerate}
\end{theorem}

\textbf{Proof. } Suppose $Z=\{u_{1}v_{1},u_1v_1,\ldots,u_{\alpha^{\prime}}%
v_{\alpha^{\prime}}\}$ is a maximum matching of $\Gamma$ and let $Y$ be the
IS of $Z$-unsaturated vertices. Observe that if $y$ and $z$ are
vertices of $Y$ and $yu_{j}\in E(\Gamma)$, then $zv_{j}\not \in E(\Gamma)$ because of
$Z$ which is maximum. Thus, for all $j\in\{1,2,\ldots,\alpha^{\prime}\}$
there are at most two edges between the sets $\{u_{j},v_{j}\}$ and $\{y,z\}$.

\begin{enumerate}
\item If there exists a vertex $v\in Y$ for which $vu_{i},vv_{i}\in E(\Gamma)$ for
some $j$, then $yu_{j},yv_{j}\not \in E(\Gamma)$ for any $y\in Y\setminus\{v\}.$
Let $A$ be the set of all vertices $v\in Y$ for which $vu_{j},vv_{j}\in E(\Gamma)$
for some $j$. Hence $xu_{j}\not \in E(\Gamma)$ or $xv_{j}\not \in E(\Gamma)$ for each
$x\in Y-A$ and any $j\in\{1,2,\ldots,\alpha^{\prime}\}$. W.l..g., we may suppose that $N(t)\subseteq\{u_{1},u_2,\ldots,u_{\alpha^{\prime
}}\}$ for any $t\in Y-A$. Define the map $l:V(\Gamma)\rightarrow\lbrack3]$ by
$l(u_{j})=2$ for $l(v_{j})=1$ for $j\in\{1,\ldots,\alpha^{\prime}\}$ and
$l(t)=1$ for any $t\in Y$. Evidently, $l$ is an OITDRDF of $\Gamma$ of weight
$p+\alpha^{\prime}(\Gamma)$ and therefore $\gamma_{tdR}^{oi}(\Gamma)\leq p+\alpha^{\prime
}(\Gamma)$.

\item Since $\Gamma$ is triangle-free, { the set} $A$ defined in item 1 is empty,
that is $vu_{j}\not \in E(\Gamma)$ or $vv_{j}\not \in E(\Gamma)$ for any $v\in Y$ and
any $j\in\{1,2,\ldots,\alpha^{\prime}\}$. Again, as before we may suppose that
$N(t)\subseteq\{u_{1},u_2,\ldots,u_{\alpha^{\prime}}\}$ for any $t\in Y$. Define
the map $l:V(\Gamma)\rightarrow\{0,1,2,3\}$ by $l(u_{j})=2$ for $l(v_{j})=1$
for $j\in\{1,2,\ldots,\alpha^{\prime}\}$ and $l(t)=0$ for any $t\in Y$. Since
$\delta(\Gamma)\geq2$, $l$ is an OITDRDF of $\Gamma$ of weight $3\alpha^{\prime}$ and
thus $\gamma_{tdR}^{oi}(\Gamma)\leq3\alpha^{\prime}(\Gamma)$.
\end{enumerate}

$\hfill\Box$

\bigskip

Since for each connected graph $\Gamma$ of order $p\ge2$, $\alpha^{\prime}(\Gamma)\le
p/2$, next result is an immediate consequence of Theorem \ref{Th1} which
improves Theorem \ref{tree1}.

\begin{corollary}
\label{3n/2}\emph{For any connected graph $\Gamma$ of order $p\geq2$, $\gamma
_{tdR}^{oi}(\Gamma)\leq3p/2$.}
\end{corollary}

\begin{rem}
\label{rmrk}Let $f=(V_{0},\ldots,V_{3})$ be an OITDRDF on a graph $\Gamma$.

\begin{enumerate}
\item As $f$ is also an OIDRDF we have that $\gamma_{oidR}(\Gamma) \le\gamma
_{tdR}^{oi}(\Gamma).$

\item Besides, every vertex $v\in V_{3}$ must have, at least, one private
neighbor in $V_{0}$ and thus $|V_{3}|\leq|V_{0}|.$
\end{enumerate}
\end{rem}

\begin{prop}
For each connected graph $\Gamma$ of order at least two,
\[
\gamma_{oidR}(\Gamma)\leq\gamma_{tdR}^{oi}(\Gamma)\leq2\gamma_{oidR}(\Gamma)-\gamma(\Gamma)
\]

\end{prop}

\textbf{Proof.} The lower bound is straightforward from Remark \ref{rmrk}. To
show the upper bound, let $l=(V_{0}^{l},\ldots,V_{3}^{l})$ be an OIDRDF of
$\Gamma$. First of all, observe that $V_{2}\cup V_{3}$ is a DS in $\Gamma$,
and so $\gamma(\Gamma)\leq|V_{2}|+|V_{3}|.$ Moreover, for every vertex $v\in
V-V_{0}^{l}$ having no neighbor in $V-V_{0}^{l}$, let $w_{v}$ be a neighbor of
$v$ belonging to $V_{0}^{l}$. Now, discuss the map $l'=(V_{0}^{l'}%
,\ldots,V_{3}^{l'})$ defined as follows $V_{2}^{l'}=V_{2}^{l}$, $V_{3}^{l'}%
=V_{3}^{l},$ $V_{1}^{l'}=V_{1}^{l}\cup\{w_{v}:v\in V-V_{0}^{l}\}$ and $l'(z)=0$
for every (eventually) not yet labeled vertex $z$. Since $\Gamma$ is connected,
$f$ is an OIDRDF, $V_{0}^{l'}\subseteq V_{0}^{l}$ and $V_{0}^{l}$ is
independent, we deduce that $l'$ is an OITDRDF on $\Gamma$. Hence
\[%
\begin{array}
[c]{rcl}%
\gamma_{tdR}^{oi}(\Gamma) & \leq & w(l')=|V_{1}^{l'}|+2|V_{2}^{l'}|+3|V_{3}%
^{l'}|\\[1em]
& \leq & \left(  |V_{1}^{l}|+|V_{2}^{l}|+|V_{3}^{l}|\right)  +2|V_{2}%
^{l}|+3|V_{3}^{l}|\\[1em]
& = & 2|V_{1}^{l}|+3|V_{2}^{l}|+4|V_{3}^{l}|\\[1em]
& \leq & 2(|V_{1}^{l}|+2|V_{2}^{l}|+3|V_{3}^{l}|)-(|V_{2}^{l}|+|V_{3}%
^{l}|)\\[1em]
& \leq & 2\gamma_{oidR}(\Gamma)-\gamma(\Gamma).
\end{array}
\]

$\hfill\Box$

\begin{prop}
Let $\Gamma$ be a connected graph of order at least two. Then
\[
\gamma_{tdR}^{oi}(\Gamma)\geq\gamma_{t,coi}(\Gamma)+\gamma(\Gamma)
\]

\end{prop}

\textbf{Proof.} Let $l=(V_{0},\ldots,V_{3})$ be an OITDRDF on $\Gamma$. Evidently,
$V_{2}\cup V_{3}$ is a DS of vertices and $V_{1}\cup V_{2}\cup
V_{3}$ is a Tco-IDS in $V(\Gamma).$ Thus,
$\gamma(\Gamma)\leq|V_{2}|+|V_{3}|$ and $\gamma_{t,coi}(\Gamma)\leq|V_{1}|+|V_{2}%
|+|V_{3}|.$ Then,
\[%
\begin{array}
[c]{rcl}%
\gamma_{t,coi}(\Gamma) & \leq & \sum_{i=1}^3|V_{i}|\\[1em]
& = & \gamma_{tdR}^{oi}(\Gamma)-\left(  |V_{2}|+2|V_{3}|\right) \\[1em]
& \leq & \gamma_{tdR}^{oi}(\Gamma)-\gamma(\Gamma)-|V_{3}|\\[1em]
& \leq & \gamma_{tdR}^{oi}(\Gamma)-\gamma(\Gamma).
\end{array}
\]
$\hfill\Box$

{ }

\begin{prop}
\label{lowerb} Let $\Gamma$ be a connected graph of order $p$. Then
\[
\gamma^{oi}_{tdR}(\Gamma) \ge\frac{2p}{\Delta}+\frac{\Delta-2}{\Delta}
\gamma_{t,coi}(\Gamma)+\frac{\gamma(\Gamma)}{\Delta}%
\]
\end{prop}

\textbf{Proof.} Suppose $h=(V_{0},V_{1},V_{2},V_{3})$ is a $\gamma_{tdR}^{oi}%
(\Gamma)$-function. As $h$ is a DRDF on $\Gamma$ we have
that any vertex in $V_{0}$ must be adjacent either to a vertex in $V_{3}$ or
either to, at least, two vertices in $V_{2}$. Let us denote by $V_{0}%
^{1}=\{v\in V_{0}:N(v)\cap V_{3}\neq\emptyset\}$ and $V_{0}^{2}=V_{0}%
\smallsetminus V_{0}^{1}.$ Since $\cup_{i=1}^3V_{i}$ is a Tco-IDS on $\Gamma$ we have that $\gamma_{t,coi}(\Gamma)\leq
\sum_{i=1}^3|V_{i}|$, $|V_{0}^{1}|\leq(\Delta-1)|V_{3}|$ and $2|V_{0}%
^{2}|\leq(\Delta-1)|V_{2}|$. So, it can be easily derived that $\Delta
|V_{2}|\geq2|V_{0}^{2}|+|V_{2}|$ and $2\Delta|V_{3}|\geq2|V_{0}^{1}%
|+2|V_{3}|.$ Moreover, as $V_{2}\cup V_{3}$ is a DS, we have that
$\gamma(\Gamma)\leq|V_{2}|+|V_{3}|$. Hence,
\[%
\begin{array}
[c]{rcl}%
\Delta\gamma_{tdR}^{oi}(\Gamma) & = & \Delta(|V_{1}|+2|V_{2}|+3|V_{3}%
|)=\Delta(\sum_{i=1}^3|V_{i}|)+\Delta|V_{2}|+2\Delta|V_{3}|\\[1em]
& \geq & \Delta(|V_{1}|+|V_{2}|+|V_{3}|)+2|V_{0}^{2}|+|V_{2}|+2|V_{0}%
^{1}|+2|V_{3}|\\[1em]
& = & (\Delta-2)(\sum_{i=1}^3|V_{i}|)+2|V_{0}|+2|V_{1}|+3|V_{2}%
|+4|V_{3}|\\[1em]
& \geq & (\Delta-2)\gamma_{t,coi}(\Gamma)+2p+|V_{2}|+2|V_{3}|\\[1em]
& \geq & (\Delta-2)\gamma_{t,coi}(\Gamma)+2p+\gamma(\Gamma)+|V_{3}|\\[1em]
& \geq & 2p+(\Delta-2)\gamma_{t,coi}(\Gamma)+\gamma(\Gamma)
\end{array}
\]
which concludes the proof. $\hfill\Box$

\bigskip

The lower bound given by Proposition \ref{lowerb} is tight, for example, for
stars and for any cycle $C_{n}$ with $n\leq10$ (except for $C_{6}$). In the
next result we give an upper bound for regular graphs that improves the known
upper bounds in some cases.

\begin{prop}
If $\Gamma$ is a connected $m$-regular graph on $p$ vertices and girth at least
$8$, then
\[
\gamma^{oi}_{tdR}(\Gamma) \le2(p-2m^{2}+3m-1)
\]

\end{prop}

\begin{figure}[ptb]
\center
\begin{tikzpicture}
\node [draw, shape=circle,fill=black,scale=0.5] (u) at (3.5,6) {};
\node at (3,6.35) {$s(u)=2$} {};
\node [draw, shape=circle,fill=black,scale=0.5] (u1) at (1,4) {};
\node at (1,4.3) {$0$} {};
\node [draw, shape=circle,fill=black,scale=0.5] (u2) at (2.5,4) {};
\node at (2.3,4.3) {$0$} {};
\node [draw, shape=circle,fill=black,scale=0.5] (ur) at (4,4) {};
\node at (3.7,4.3) {$0$} {};
\node [draw, shape=circle,fill=black,scale=0.5] (v) at (7,6) {};
\node at (8,6.35) {$s(v)=2$} {};
\node [draw, shape=circle,fill=black,scale=0.5] (v1) at (6,4) {};
\node at (5.85,4.3) {$0$} {};
\node [draw, shape=circle,fill=black,scale=0.5] (v2) at (7.5,4) {};
\node at (7.2,4.3) {$0$} {};
\node [draw, shape=circle,fill=black,scale=0.5] (vr) at (9,4) {};
\node at (9,4.3) {$0$} {};
\node [draw, shape=circle,fill=black,scale=0.5] (z11) at (.5,2) {};
\node at (.25,2.25) {$2$} {};
\node [draw, shape=circle,fill=black,scale=0.5] (z12) at (1,2) {};
\node at (.85,2.25) {$1$} {};
\node [draw, shape=circle,fill=black,scale=0.5] (z1r) at (1.5,2) {};
\node at (1.35,2.25) {$1$} {};
\node [draw, shape=circle,fill=black,scale=0.5] (z21) at (2,2) {};
\node at (1.85,2.25) {$2$} {};
\node [draw, shape=circle,fill=black,scale=0.5] (z22) at (2.5,2) {};
\node at (2.35,2.25) {$1$} {};
\node [draw, shape=circle,fill=black,scale=0.5] (z2r) at (3,2) {};
\node at (2.85,2.25) {$1$} {};
\node [draw, shape=circle,fill=black,scale=0.5] (zr1) at (3.5,2) {};
\node at (3.35,2.25) {$2$} {};
\node [draw, shape=circle,fill=black,scale=0.5] (zr2) at (4,2) {};
\node at (3.85,2.25) {$1$} {};
\node [draw, shape=circle,fill=black,scale=0.5] (zrr) at (4.5,2) {};
\node at (4.35,2.25) {$1$} {};
\node [draw, shape=circle,fill=black,scale=0.5] (t11) at (5.5,2) {};
\node at (5.25,2.25) {$2$} {};
\node [draw, shape=circle,fill=black,scale=0.5] (t12) at (6,2) {};
\node at (5.85,2.25) {$1$} {};
\node [draw, shape=circle,fill=black,scale=0.5] (t1r) at (6.5,2) {};
\node at (6.35,2.25) {$1$} {};
\node [draw, shape=circle,fill=black,scale=0.5] (t21) at (7,2) {};
\node at (6.85,2.25) {$2$} {};
\node [draw, shape=circle,fill=black,scale=0.5] (t22) at (7.5,2) {};
\node at (7.35,2.25) {$1$} {};
\node [draw, shape=circle,fill=black,scale=0.5] (t2r) at (8,2) {};
\node at (7.85,2.25) {$1$} {};
\node [draw, shape=circle,fill=black,scale=0.5] (tr1) at (8.5,2) {};
\node at (8.35,2.25) {$2$} {};
\node [draw, shape=circle,fill=black,scale=0.5] (tr2) at (9,2) {};
\node at (8.85,2.25) {$1$} {};
\node [draw, shape=circle,fill=black,scale=0.5] (trr) at (9.5,2) {};
\node at (9.35,2.25) {$1$} {};
\node [draw, shape=circle,fill=black,scale=0.5] (a11) at (0.25,0) {};
\node at (.25,-0.35) {$1$} {};
\node [draw, shape=circle,fill=black,scale=0.5] (a12) at (0.5,0) {};
\node at (.5,-0.35) {$1$} {};
\node [draw, shape=circle,fill=black,scale=0.5] (a1r) at (0.75,0) {};
\node at (0.75,-0.35) {$1$} {};
\node [draw, shape=circle,fill=black,scale=0.5] (a21) at (1.75,0) {};
\node at (1.75,-0.35) {$1$} {};
\node [draw, shape=circle,fill=black,scale=0.5] (a22) at (2,0) {};
\node at (2,-0.35) {$1$} {};
\node [draw, shape=circle,fill=black,scale=0.5] (a2r) at (2.25,0) {};
\node at (2.25,-0.35) {$1$} {};
\node [draw, shape=circle,fill=black,scale=0.5] (ar1) at (3.25,0) {};
\node at (3.25,-0.35) {$1$} {};
\node [draw, shape=circle,fill=black,scale=0.5] (ar2) at (3.5,0) {};
\node at (3.5,-0.35) {$1$} {};
\node [draw, shape=circle,fill=black,scale=0.5] (arr) at (3.75,0) {};
\node at (3.75,-0.35) {$1$} {};
\node [draw, shape=circle,fill=black,scale=0.5] (b11) at (5.25, 0) {};
\node at (5.25,-0.35) {$1$} {};
\node [draw, shape=circle,fill=black,scale=0.5] (b12) at (5.5,0) {};
\node at (5.5,-0.35) {$1$} {};
\node [draw, shape=circle,fill=black,scale=0.5] (b1r) at (5.75,0) {};
\node at (5.75,-0.35) {$1$} {};
\node [draw, shape=circle,fill=black,scale=0.5] (b21) at (6.75,0) {};
\node at (6.75,-0.35) {$1$} {};
\node [draw, shape=circle,fill=black,scale=0.5] (b22) at (7,0) {};
\node at (7,-0.35) {$1$} {};
\node [draw, shape=circle,fill=black,scale=0.5] (b2r) at (7.25,0) {};
\node at (7.25,-0.35) {$1$} {};
\node [draw, shape=circle,fill=black,scale=0.5] (br1) at (8.25,0) {};
\node at (8.25,-0.35) {$1$} {};
\node [draw, shape=circle,fill=black,scale=0.5] (br2) at (8.5,0) {};
\node at (8.5,-0.35) {$1$} {};
\node [draw, shape=circle,fill=black,scale=0.5] (brr) at (8.75,0) {};
\node at (8.75,-0.35) {$1$} {};
\draw (z11)--(u1)--(u)--(v)--(v1)--(t11);
\draw (u)--(u2)--(z21);
\draw (u)--(ur)--(zr1);
\draw (z12)--(u1)--(z1r);
\draw (z22)--(u2)--(z2r);
\draw (zr2)--(ur)--(zrr);
\draw (v)--(v2)--(t21);
\draw (v)--(vr)--(tr1);
\draw (t12)--(v1)--(t1r);
\draw (t22)--(v2)--(t2r);
\draw (tr2)--(vr)--(trr);
\draw (z11)--(a11)--(z11)--(a12)--(z11)--(a1r);
\draw (z21)--(a21)--(z21)--(a22)--(z21)--(a2r);
\draw (zr1)--(ar1)--(zr1)--(ar2)--(zr1)--(arr);
\draw (t11)--(b11)--(t11)--(b12)--(t11)--(b1r);
\draw (t21)--(b21)--(t21)--(b22)--(t21)--(b2r);
\draw (tr1)--(br1)--(tr1)--(br2)--(tr1)--(brr);
\draw (5,-2) ellipse (5cm and 1 cm );
\node at (5,-2) {$s(x)=2, \text{ otherwise. }$ } {};
\end{tikzpicture}
\caption{An OITDRD function in an $r$-regular graph.}%
\label{regular}%
\end{figure}
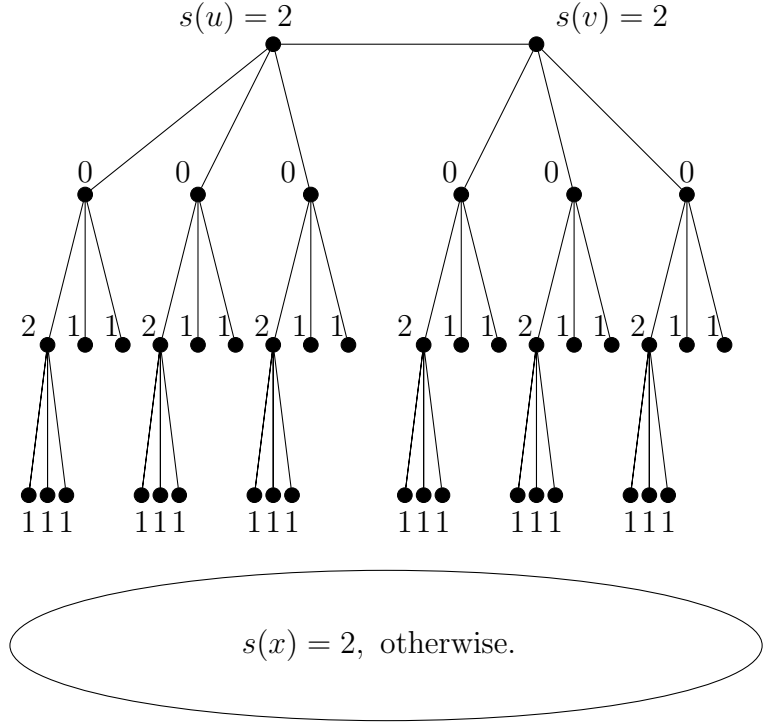

\textbf{Proof.} Suppose $r$ and $r'$ are two adjacent vertices in $V(\Gamma)$, and let us denote
by $N(r')=\{r,r'_{1},r'_2,\ldots, r'_{m-1} \}, N(r)=\{ r',r_{1},r_2,\ldots,r_{m-1} \},
N(r'_{j})=\{ r',z_{j}^{1},z_j^2,\ldots,z_{j}^{m-1} \}$ and finally $N(r_{j})=\{
r,t_{j}^{1},t_j^2,\ldots,t_{j}^{m-1} \}.$ With these notations, let $s$ be the
map defined as $s(r'_{j})=s(r_{j})=0, s(z_{k})=s(t_{k})=s(x)=1$ for $1\le
j \le m-1, 2\le k\le m-1$ and for all $x\in(N(z_{j}^{1})-r'_{j}) \cup
(N(t_{j}^{1})-r_{j})$, and $s(x)=2$ otherwise (see Figure~\ref{regular}).
Since the girth is at least $8$, these sets of vertices are disjoint. Every
vertex labeled with a $0$ is adjacent to two vertices labeled with $2$ and
each vertex with a label $1$ has a neighbor with a label $2$, therefore $s$ is
a DRDF on $\Gamma$. Evidently the vertices with a label $0$ induce an edgeless
subgraph and no vertex with positive label is isolated. Hence, $s$ is an
OITDRDF on $\Gamma$ and
\[%
\begin{array}
[c]{rcl}%
\gamma^{oi}_{tdR}(\Gamma) & \le & \displaystyle s(r')+s(r)+ \sum_{
\begin{array}
[c]{c}%
\scriptstyle x\in N(r')\cup N(r)\\[-.4em]%
\scriptstyle x \neq r,r'
\end{array}
} s(x) + \sum_{%
\begin{array}
[c]{c}%
\scriptstyle x\in N(z_{j}^{1})-r'_{j}\\[-.4em]%
\scriptstyle x\in N(t_{j}^{1})-r_{j}%
\end{array}
} s(x) + \sum_{otherwise} s(x)\\[0.5 em]
& = & 4+2m(m-1)+2(m-1)^{2}+2(p-2-2(m-1)-2(m-1)^{2}-2(m-1)^{2})\\[0.5 em]
& = & 2(p-2m^{2}+3m-1)
\end{array}
\]

This bound is tight for $C_{8}$ and improves the bound given by Corollary
\ref{3n/2} for $m$-regular graphs whenever $p< 4(2m^{2}-3m+1).$

\section{Trees}

In this part we present tight bounds on the OITDRD-number of trees. First
we characterize the trees attaining the bound of Corollary \ref{3n/2}. Suppose
that $\mathcal{M}=\{Cor(M)\mid M \; \text{is a tree} \}$. We start with a
simple result.

\begin{corollary}
\label{c1} \emph{For $t\ge r\ge1$, $\gamma_{oitdR}(DS_{r,t})=6$.}
\end{corollary}

\textbf{Proof. } Let $z_{1},z_{2}$ be the non-leaf vertices of $DS_{r,t}$ and
$w_{i}$ a leaf neighbor of $z_{j}$ for any $j\in\{1,2\}$. By Observation \ref{ob},
we find $h(w_{j})+h(z_{j})\geq3$ for any $j\in\{1,2\}$ and so $\gamma_{tdR}%
^{oi}(DS_{r,t})\geq6$. On the other hand, the map $h=(V(DS_{r,t}%
)-\{z_{1},z_{2}\},\emptyset,\emptyset,\{z_{1},z_{2}\})$ is evidently a OITDRDF
of $DS_{r,t}$ yielding $\gamma_{tdR}^{oi}(DS_{r,t})\leq6$. Thus $\gamma
_{tdR}^{oi}(DS_{r,t})=6$. $\hfill\Box$

We apply ``In-Hyp'' instead of ``Induction hypothesis''.

\begin{theorem}
\label{tree}\emph{Let $M$ be a tree of order $p\geq2$. Then $\gamma_{tdR}%
^{oi}(M)=\frac{3p}{2}$ iff $T\in\mathcal{M}$.}
\end{theorem}

\textbf{Proof. } If $M\in\mathcal{T}$, then by Observation \ref{ob},
$\gamma_{tdR}^{oi}(M)\geq\frac{3p}{2}$ and the equality follows from Corollary
\ref{3n/2}. Suppose now that $\gamma_{tdR}^{oi}(M)=\frac{3p}{2}$. To show that
$M\in\mathcal{T}$, we use an induction on $p$. The result is trivial for
$p\leq3$. Assume $p\geq4$ and the result holds for all trees of order at most $p-1$.
Let $M$ be a tree of order $p$. If $diam(M)=2$, then
$M$ is a star $K_{1,p-1}$ and evidently $\gamma_{tdR}^{oi}(K_{1,p-1}%
)=4<\frac{3p}{2}$. If $diam(M)=2$, then $M$ is $DS_{r,t}$ for some $t\geq r\geq1$ and by Corollary \ref{c1}, we have
$\gamma_{tdR}^{oi}(DS_{r,t})=6\leq\frac{3(r+t+2)}{2},$ with equality iff
$r=t=1$, which is $M=P_{4}\in\mathcal{T}$. Therefore we suppose the diameter of $M$
is greater than 3. Let $P:=z_{1}z_{2}\ldots z_{d}\;(d\geq5)$ be a diametrical
path on $M$ and root $M$ at $z_{k}$. Assume $M^{\prime}=T-T_{z_{2}}$. By
Corollary \ref{3n/2}, we have $\gamma_{tdR}(M^{\prime})\leq\frac{3p(M^{\prime
})}{2}$. If $d_{T}(z_{2})\geq3$, then each $\gamma_{tdR}^{oi}(M^{\prime}%
)$-function $h$ that may extend to an OITDRDF of $M$ by labeling a 3 to
$z_{2}$, a 1 to $z_{1}$ and the other leaf neighbors of $z_{2}$ label by 0, which we find
\[
\begin{array}{clc}
\gamma_{tdR}^{oi}(M)& \leq & \gamma_{tdR}^{oi}(M^{\prime})+4\hfill\\
& \leq& \frac{3(p-d_{T}(z_{2}))}{2}+4\hfill\\
& < & \frac{3p}{2}\hfill\\.
\end{array}
\]
Hence we can assume that $d_{T}(z_{2})=2$. Then each $\gamma_{tdR}%
^{oi}(M^{\prime})$-function $h$ that may extend to an OITDRDF of $M$ by
labeling a 2 to $z_{2}$ and a 1 to $z_{1}$, and we find
\begin{equation}
\gamma_{tdR}^{oi}(M)\leq\gamma_{tdR}^{oi}(M^{\prime})+3\leq\frac{3(p-2)}%
{2}+3=\frac{3p}{2}. \label{eqp}%
\end{equation}
The equality holds in (\ref{eqp}), iff $\gamma_{tdR}^{oi}(M^{\prime}%
)=\frac{3p(M^{\prime})}{2}$. By the In-Hyp we find $M^{\prime}%
\in\mathcal{T}$. Suppose that $M^{\prime}=Cor(M_{1})$. We claim that $z_{3}\in
V(M_{1})$. Assume, to the contrary, that $z_{3}\not \in V(M_{1})$. Then
$z_{3}$ is a leaf of $M^{\prime}$. Define $t:V(\Gamma)\rightarrow [3]$ by
$t(z)=2$ for $z\in V(M_{1})$ and $t(z)=1$ otherwise. Evidently $t$ is a
$\gamma_{tdR}(M^{\prime})$-function. Now the map $t':V(M)\rightarrow
[3]$ defined by $t'(z_{3})=0,t'(z_{1})=1,t'(z_{2})=2$ and $t'(z)=t(z)$
otherwise, is evidently is an OITDRDF of $M$ of weight less than $\frac{3p}{2}$
yielding a contradiction. Thus $z_{3}\in V(M_{1})$. Suppose $M_{2}$ is the tree
obtained from $M_{1}$ by adding the vertex $z_{2}$ and the edge $z_{3}z_{2}.$
In this case, one can see that $M=Cor(M_{2})\in\mathcal{T}$ and this finishes the prove. $\hfill\Box$

\bigskip

In the aim to state our next upper bound, we make the following two
observations whose proofs are straightforward.

\begin{obs}
\label{obs1} Let $M$ be a tree with diameter greater than or equal to 3 and let $u$ be a
strong stem of $M$. Then there exists a $\gamma_{tdR}^{oi}(M)$-function $h$
for which $h(u)=3$ and $h(z)=0$ for every leaf adjacent to $u$.
\end{obs}

\begin{obs}
\label{obs2} Let $M$ be a tree and let $r$ be a weak stem { with a }unique
leaf $r'$. Then there exists a $\gamma_{tdR}^{oi}(M)$-function $t$ for which
$t(r)+t(r')=3$.
\end{obs}

\begin{theorem}
\label{thm3}For every tree $M$ of order $p\geq3$ with $s(M)$ stems,
\[
\gamma_{tdR}^{oi}(M)\leq\frac{6p(M)+3s(M)}{5}.
\]
This bound is tight for each tree $M\in\mathcal{T}$.
\end{theorem}

\textbf{Proof. } We proceed by induction on $p$. If $p=3$, then $M=P_{3}$
and $\gamma_{tdR}^{oi}(M)=4<\frac{6p(M)+3s(M)}{5}$. If $p=4$, then
either $M=P_{4}$, then $M\in\mathcal{T}$ and $\gamma_{tdR}^{oi}%
(M)=6=\frac{6p(M)+3s(M)}{5}$  or $M=K_{1,3}$  and $\gamma_{tdR}%
^{oi}(M)=4<\frac{6p(M)+3s(M)}{5}$. Hence let $p\geq5$ and suppose that
for any tree $M^{\prime}$ of order $p(M^{\prime})$ with $s(M^{\prime})$ stems
and $3\leq p(M^{\prime})<p$, { satisfies }$\gamma_{tdR}^{oi}(M^{\prime}%
)\leq\frac{6p(M^{\prime})+3s(M^{\prime})}{5}$. Suppose $M$ is a tree of order
$p\geq5$ with $s(M)$ stems. If $\mathrm{diam}(M)=2$, { then }$M=K_{1,p-1}$ {
and }$\gamma_{tdR}^{oi}(M)=4<\frac{6p(M)+3s(M)}{5}$. If $\mathrm{diam}(M)=3$,
 then $M$ is a double star  and $\gamma_{tdR}^{oi}(M)=6<\frac
{6p(M)+3s(M)}{5}$.  Hence in the sequel, we may suppose that $\mathrm{diam}%
(M)\geq4$.

First let $M$ has a strong stem $r$,  and suppose $w$ { is a leaf neighbor }of
$r$.  Suppose $M^{\prime}$  is the tree obtained from $M$ { by deleting }$w$
and let $t$ be a $\gamma_{tdR}^{oi}(M^{\prime})$-function. Note that $r$
remains a stem of $M^{\prime}$. If $r$ is a weak stem of $M^{\prime}$  with
a unique leaf $r'$, then by Observation \ref{obs2}, we may suppose that
$t(r)+t(r')=3,$  with in addition, without loss of generality, $t(r')<t(r)$.
Moreover, if $r$ is a strong stem of $M^{\prime}$, then by Observation
\ref{obs1}, we may suppose that $t(r)=3$ and $t(r'')=0$ for any leaf  neighbor
of $r$ in $M^{\prime}$. Hence in either case, the map $t'$ defined {
on }$V(M)$ by $t'(w)=1$ and $t'(r)=t(r)$ for each $r\in V(M^{\prime})$, is an
OITDRDF on $M$. Applying the induction on $M^{\prime}$ and applying the facts
that $p=p(M^{\prime})+1$ and $s(M)=s(M^{\prime})$, we find
\begin{align}
\gamma_{tdR}^{oi}(M)  &  \leq t'(V(M))\nonumber\\
&  =t(V(M^{\prime}))+t(w)\nonumber\\
&  \leq\frac{6p(M^{\prime})+3s(M^{\prime})}{5}+1\nonumber\\
&  =\frac{6(p-1)+3s(M)}{5}+1\nonumber\\
&  <\frac{6p+3s(M)}{5}.\nonumber
\end{align}
 Therefore, we can suppose that each stem of $M$ is weak. If $M$ is a path,
then by Lemma \ref{pathandcycle}, we have $\gamma_{tdR}^{oi}%
(M)=\left\lceil \frac{6p(M)}{5}\right\rceil <\frac{6p(M)+3s(M)}{5},$ as
desired. Hence we may suppose that $M$ is not a path, that is, there exists
some vertex of degree at least $3$ in $M$. Let $P=r_{1}r_{2}\dots r_{d+1}$ be
a diametral path of $M$, where $d=\mathrm{diam}(M)\geq4$, such that $r_{k}$ is
the first vertex of $P$ with $d_{T}(r_{k})\geq3$. By our earlier assumption
that every stem is weak, we deduce that $k\geq3$. We discuss the following situations.

\smallskip\noindent\textbf{Case 1.} Assume that $k=3$.\newline Suppose $M^{\prime
}=T-\{r_{1},r_{2}\}$ and let $t$ be a $\gamma_{tdR}^{oi}(M^{\prime}%
)$-function. Then the map $t'$ defined { on }$V(M)$ by $t'(r_{1})=1$,
$t'(r_{2})=2$ and $t'(r)=t(r)$ for any $r\in V(M^{\prime})$, is an OITDRDF on
$M$. Observe that $p=p(M^{\prime})+2$ and $s(M)=s(M^{\prime})+1$. By the
In-Hyp on $M^{\prime}$, we find
\begin{align}
\gamma_{tdR}^{oi}(M)  &  \leq t'(V(M))\nonumber\\
&  =t(V(M^{\prime}))+t'(r_{1})+t'(r_{2})\nonumber\\
&  \leq\frac{6p(M^{\prime})+3s(M^{\prime})}{5}+3\nonumber\\
&  =\frac{6(p-2)+3(s(M)-1)}{5}+3\nonumber\\
&  =\frac{6p+3s(M)}{5}.\nonumber
\end{align}

\smallskip\noindent\textbf{Case 2.} $k=4$.\newline Let $M^{\prime}%
=T-\{r_{1},r_{2},r_{3}\}$ and let $t$ be a $\gamma_{tdR}^{oi}(M^{\prime}%
)$-function. Then the map $t'$ defined { on }$V(M)$ by $t'(r_{1}%
)=t'(r_{3})=1$, $t'(r_{2})=2$ and $t'(r)=t(r)$ for any $r\in V(M^{\prime})$, is
an OITDRDF on $M$. Observe that $p=p(M^{\prime})+3$ and $s(M)=s(M^{\prime}%
)+1$. By the In-Hyp, we find
\begin{align}
\gamma_{oitdR}(M)  &  \leq t'(V(M))\nonumber\\
&  =t(V(M^{\prime}))+\sum_{i=1}^3t'(r_{i})\nonumber\\
&  \leq\frac{6p(M^{\prime})+3s(M^{\prime})}{5}+4\nonumber\\
&  =\frac{6(p-3)+3(s(M)-1)}{5}+4\nonumber\\
&  <\frac{6p+3s(M)}{5}.\nonumber
\end{align}

\smallskip\noindent\textbf{Case 3.} $k=5$.\newline First, suppose that $r_{5}$
is a stem in $M$. Suppose $M^{\prime}$ { is the tree obtained from }$M$ { by
deleting vertices }$r_{1},r_{2},r_{3},r_{4},$ and let $t$ be a $\gamma
_{oitdR}(M^{\prime})$-function.  Evidently, $r_{5}$  is a stem in
$M^{\prime}$.  Suppose $w$ { is the unique leaf neighbor of }$r_{5}$.  Note
that we may suppose that $w$  does not belong to $P$  since $d_{T}%
(r_{5})\geq3.$ By Observation \ref{obs2}, { we assume that }$t(r_{5}%
)+t(w)=3$, with in addition $t(r_{5})\geq t(w)$. Then the map $t'$
defined  on $V(M)$ by $t'(r_{1})=1$, $t'(r_{2})=t'(r_{3})=2$, $t'(r_{4})=0$ and
$t'(r)=t(r)$ for any $r\in V(M^{\prime})$, is an OITDRDF on $M$. Observe that
$p=p(M^{\prime})+4$ and $s(M)=s(M^{\prime})+1$. By the In-Hyp, we find
\begin{align}
\gamma_{tdR}^{oi}(M)  &  \leq t'(V(M))\nonumber\\
&  =t(V(M^{\prime}))+\sum_{i=1}^4t'(r_{i})\nonumber\\
&  \leq\frac{6p(M^{\prime})+3s(M^{\prime})}{5}+5\nonumber\\
&  =\frac{6(p-4)+3(s(M)-1)}{5}+5\nonumber\\
&  <\frac{6p+3s(M)}{5}.\nonumber
\end{align}

Second, assume that $r_{5}$ is not a stem in $M$. Assume $M_{1}$ is the
component of $M-r_{6}$ containing $r_{5}$. If $M_{1}$  contains a vertex of
degree at least three in $M_{1}$ other than $r_{5},$  then by a similar
method to that applied in Case 2 or 3, we get $\gamma_{tdR}^{oi}%
(M)<\frac{6p+3s(M)}{5}$. Hence we may assume that  every vertex of
$M_{1}$  but possibly $r_{5}$  has degree $1$ or $2$. If there exists some
leaf at distance $2$ or $3$ from $r_{5}$ in $M_{1}$, then again by a similar
 method to that applied in Case 2 or 3, we obtain $\gamma_{tdR}^{oi}%
(M)<\frac{6p+3s(M)}{5}$.  Therefore, in what follows we can suppose that
each leaf in $M_{1}$  is at distance four from $r_{5},$  that is, %
$M_{1}$ { is a tree obtained from a star }$K_{1,t-1}$ of order $t\geq3$ { and
centered at }$r_{5}$ by { subdividing any edge of the star }exactly three
times. { Let }$M_{2}=T-T_{1}$ and let $t$ be a $\gamma_{tdR}^{oi}(M_{2}%
)$-function. { Note that since }$P=r_{1}r_{2}\dots r_{d+1}$ is a diametral
path, $p(M_{2})\geq4$.  We also note that the map $t'$ defined by
$t'(r)=2$ for every vertex $r$ at distance $1$ or $3$ from $r_{5}$ in $M_{1}$,
$t'(r)=0$ for every vertex $r$ at distance $2$ from $r_{5}$ in $M_{1}$,
$t'(r)=1$ for { the }remaining vertices $r$ of $M_{1}$ and $t'(r)=t(r)$ for
$r\in V(M_{2})$, is an OITDRDF on $M$. Observe that $p=p(M_{2})+4t-3$ and
$s(M_{2})+t-2\leq s(M)\leq s(M_{2})+t-1$. By the In-Hyp { on }$M_{2}$, we
obtain
\begin{align}
\gamma_{tdR}^{oi}(M)  &  \leq t'(V(M))\nonumber\\
&  =t'(V(M_{1}))+t(V(M_{2}))\nonumber\\
&  \leq5(t-1)+1+\frac{6p(M_{2})+3s(M_{2})}{5}\nonumber\\
&  \leq5t-4+\frac{6(p-4t+3)+3(s(M)-t+2)}{5}\nonumber\\
&  <\frac{6p+3s(M)}{5}.\nonumber
\end{align}

\smallskip\noindent\textbf{Case 4.} $k\geq6$.\newline Let $M^{\prime
}=T-\{r_{1},r_{2},r_{3},r_{4},r_{5}\}$ and let $t$ be a $\gamma_{tdR}%
^{oi}(M^{\prime})$-function. Then the map $t'$ defined  on $V(M)$ by
$t'(r_{1})=t'(r_{5})=1$, $t'(r_{2})=t'(r_{4})=2$, $t'(r_{3})=0$ and $t'(r)=t(r)$ for
each $r\in V(M^{\prime})$, is an OITDRDF on $M$. Observe that $p=p(M^{\prime
})+5$ and $s(M^{\prime})\leq s(M)$. By the In-Hyp { on }$M^{\prime}$, we find
\begin{align}
\gamma_{tdR}^{oi}(M)  &  \leq t'(V(M))\nonumber\\
&  =t(V(M^{\prime}))+\sum_{i=1}^{5}t'(r_{i})\nonumber\\
&  \leq\frac{6p(M^{\prime})+3s(M^{\prime})}{5}+6\nonumber\\
&  =\frac{6(p-5)+3s(M)}{5}+6\nonumber\\
&  \leq\frac{6p+3s(M)}{5},\nonumber
\end{align}
as desired. This finishes the prove. $\hfill\Box$

\end{document}